\documentclass[11pt]{amsart}

\usepackage[english]{babel}
\usepackage[OT2,T1]{fontenc}
\usepackage[utf8]{inputenc}

\usepackage{mathpazo}
\usepackage{graphicx}
\usepackage{microtype}
\usepackage{enumerate}
\usepackage{fontawesome}

\usepackage{csquotes}
\usepackage[dvipsnames]{xcolor}
\usepackage{url}

\usepackage[breaklinks,bookmarks, colorlinks, citecolor=MidnightBlue!75!black, urlcolor=MidnightBlue!75!black, linkcolor=MidnightBlue!75!black]{hyperref}

\calclayout

\AtBeginDocument{%
   \def\MR#1{}
}

\usepackage{scalerel}

\usepackage{amsmath, amssymb, amsthm, amscd}
\usepackage{commath, mathtools, mathrsfs}
\usepackage{thmtools}
\usepackage{tikz, tikz-cd}
\usetikzlibrary{matrix, arrows, decorations.pathmorphing}

\theoremstyle{plain}

\newtheorem{theorem}{Theorem}[section]
\newtheorem{lemma}[theorem]{Lemma}
\newtheorem{proposition}[theorem]{Proposition}
\newtheorem{corollary}[theorem]{Corollary}

\newtheorem*{question*}{Question}

\newtheorem*{theorem*}{Theorem}
\newtheorem*{lemma*}{Lemma}
\newtheorem*{proposition*}{Proposition}
\newtheorem*{corollary*}{Corollary}

\declaretheorem[numbered=no,name=Theorem A]{TA}
\declaretheorem[numbered=no,name=Theorem B]{TB}

\declaretheorem[numbered=no,name=Corollary A]{CA}
\declaretheorem[numbered=no,name=Corollary B]{CB}
\declaretheorem[numbered=no,name=Corollary C]{CC}

\theoremstyle{definition}

\newtheorem{definition}[theorem]{Definition}
\newtheorem{example}[theorem]{Example}

\newtheorem*{conjecture*}{Conjecture}
\newtheorem*{remark*}{Remark}
\newtheorem*{definition*}{Definition}
\newtheorem*{observation*}{Observation}

\newcommand{\on}{\operatorname}

\newcommand{\A}{\mathbb{A}}

\newcommand{\aut}{\on{Aut}}

\newcommand{\C}{\mathbb{C}}

\newcommand{\gal}{\on{Gal}}

\newcommand{\hz}{\hat{\mathbb{Z}}}

\newcommand{\mc}{\mathcal}

\newcommand{\Q}{\mathbb{Q}}

\newcommand{\s}{\subseteq}

\newcommand{\spec}{\on{Spec}}

\newcommand{\xar}{\xrightarrow}
\newcommand{\Z}{\mathbb{Z}}

\renewcommand{\epsilon}{\varepsilon}
\renewcommand{\H}{\on{H}}
\renewcommand{\hom}{\on{Hom}}

\renewcommand{\P}{\mathbb{P}}
\renewcommand{\phi}{\varphi}
\renewcommand{\projlim}{\varprojlim}

\renewcommand{\tilde}{\widetilde}
\renewcommand{\hat}{\widehat}

\DeclareSymbolFont{cyrletters}{OT2}{wncyr}{m}{n}
\DeclareMathSymbol{\Sha}{\mathalpha}{cyrletters}{"58}

\usepackage[giveninits=true,style=alphabetic,backend=biber,maxnames=99,maxalphanames=5,giveninits=true]{biblatex}
\bibliography{scfg}

\usepackage{todonotes}

\declaretheorem[numbered=no,name=Weak Bombieri-Lang Conjecture]{BL}
\declaretheorem[numbered=no,name=Geometric Lang Conjecture]{GL}

\usepackage{graphicx}
\usepackage{microtype}

\author{Giulio Bresciani}
\address[G. Bresciani]{CRM Ennio de Giorgi, Collegio Puteano, Office 21, Piazza dei Cavalieri 3, 56126 Pisa}
\email{giulio.bresciani@gmail.com}
\thanks{The author was partially supported by the DFG-funded Priority Program "Homotopy Theory and Algebraic Geometry" SPP 1786}

\title{On the section conjecture over fields of finite type}
\date{}

\begin{document}

\begin{abstract}
	Assume that the section conjecture holds over number fields. We prove then that it holds for a broad class of curves defined over finitely generated extensions of $\mathbb{Q}$. This class contains every projective, hyperelliptic curve, every hyperbolic, affine curve of genus $\le 2$, and a non-empty open subset of any curve. If we furthermore assume the weak Bombieri-Lang conjecture, we prove that the section conjecture holds for every hyperbolic curve over every finitely generated extension of $\mathbb{Q}$. 
\end{abstract}

\maketitle

\section{Introduction}

Curves are assumed to be smooth and geometrically connected, but not necessarily projective. Fields are assumed to be of characteristic $0$.

\subsection{Grothendieck's section conjecture} Let $X$ be a hyperbolic curve over a field $k$ finitely generated over $\Q$ with completion $\bar{X}\supseteq X$, there is a short exact sequence of étale fundamental groups
\[0\to\pi_{1}(X_{\bar{k}})\to\pi_{1}(X)\to\gal(\bar{k}/k)\to 0.\]
The set of sections of this sequence modulo the action of $\pi_{1}(X_{\bar{k}})$ is called the \emph{space of Galois sections} of $X$. A rational point of $X$ induces a Galois section, while a rational point in the boundary $\bar{X}\setminus X$ induces an infinite number of them \cite{eh08}. A Galois section is \emph{geometric} if it is associated to a point of $X(k)$, \emph{cuspidal} if it is associated to a point of $\bar{X}\setminus X(k)$, and it \emph{comes from geometry} if it is either geometric or cuspidal.

Grothendieck's section conjecture \cite{gro97} states that every Galois section of $X$ comes from geometry. Most of the known results about the conjecture are for number fields and their completions \cite{koe05} \cite{sti10} \cite{sti15} \cite{cs15} \cite{bre22}. The following is the subject of this paper.

\begin{question*}
	Assume that the section conjecture holds for all hyperbolic curves over all number fields. Does it hold for all hyperbolic curves over all fields finitely generated over $\Q$?
\end{question*}

We always assume that the section conjecture holds over number fields for some class of hyperbolic curves. When we say that some result is "conditional" or "unconditional" we are referring to some other conjecture, since the section conjecture over number fields is always given for granted.

\subsection{Known results}

M. Sa\"idi proved two results concerning the question above in \cite{sai16}. The first is conditional on a conjecture similar to the Tate-Shafarevich conjecture and regards all hyperbolic curves, while the second is unconditional but restricted to curves which descend to number fields. Furthermore, in joint paper with M. Tyler \cite{saty21}, he gave a positive answer to the analogous question for a birational version of the section conjecture \cite[\S 18.4]{sti13}.

If $C$ is a curve over a field $k$ and $A$ is an abelian variety over $K=k(C)$, Sa\"idi defines the Tate-Shafarevich group $\Sha(C,A)$ as the kernel of $\H^{1}(K,A)\to\prod_{c\in C^{(1)}}\H^{1}(K_{c},A_{K_{c}})$ analogously to the case of number fields. Let $\Sha(C,A)[n]$ be the $n$-torsion subgroup, and write $T\Sha(C,A)=\projlim_{n}\Sha(C,A)[n]$ for the Tate module.

\begin{theorem*}[{\cite[Corollary B]{sai16}}]
	Assume that the section conjecture holds for all projective, hyperbolic curves over all number fields, and that $T\Sha(C,A)=0$ for every $C,A,k$ as above with $k$ finitely generated over $\Q$. Then the section conjecture holds for all projective, hyperbolic curves over all fields finitely generated over $\Q$.
\end{theorem*} 

\begin{theorem*}[{\cite[Theorem C]{sai16}}]
	Assume that the section conjecture holds for all projective, hyperbolic curves over all number fields. Let $K$ be a finitely generated extension of $\Q$ with algebraic closure $\bar{K}$, and $\bar{\Q}$ the algebraic closure of $\Q$ in $\bar{K}$. Then the section conjecture holds for every projective, hyperbolic curve $X$ over $K$ such that $X_{\bar{K}}$ descends to $\bar{\Q}$.
\end{theorem*} 

Sa\"idi's work \cite{sai16} is centered on the study of Tate-Shafarevich group of the relative Jacobian of a relative curve. We study directly the family of curves without passing through the relative Jacobian. We obtain a different conditional result and strengthen the unconditional result. In particular, our unconditional result is strong enough to give a complete answer in the case of hyperelliptic curves and affine curves of genus $\le 2$.

\subsection{Specialization of curves}

In order to better state our results, let us define what it means for a curve to be a \emph{specialization} of another curve. Let $X_{1},X_{2}$ be curves of genus $g$ over fields $k_{1},k_{2}$. We say that they are \emph{equivalent} if there exists a field $K$ containing copies of $k_{1},k_{2}$ such that $X_{1,K}\simeq X_{2,K}$.

We say that $X_{2}$ is a \emph{specialization} of $X_{1}$ if there exists an integral scheme $S$, a smooth, projective family of curves $\bar{\mc{X}}\to S$, a divisor $D\subset \bar{\mc{X}}$ finite étale over $S$ such that $X_{1}$ is equivalent to the generic fiber of $\bar{\mc{X}}\setminus D$ and $X_{2}$ is equivalent to some fiber.

Equivalently, up to enlarging $k_{1},k_{2}$, there exist curves $\bar{X}_{i}$ over $k_{i}$ and points $x_{i,1},\dots,x_{i,n}\in\bar{X}_{i}$ with $X_{i}=\bar{X}_{i}\setminus\{x_{i,1},\dots,x_{i,n}\}$ such that the marked curve $(\bar{X}_{2},x_{2,1},\dots,x_{2,n})$ is in the scheme theoretic closure of $(\bar{X}_{1},x_{1,1},\dots,x_{1,n})$ in the moduli stack of marked curves $\mc{M}_{g,n}/\Q$. Specialization defines a partial order on the class of curves modulo equivalence.


\subsection{New, conditional result}


A geometrically connected variety is of general type if it is birational to a smooth, projective variety whose Kodaira dimension is equal to its Krull dimension (the Kodaira dimension is a birational invariant). In dimension $1$, curves of general type are those of genus $\ge 2$. Let us recall two conjectures of Bombieri and Lang.

\begin{BL}
	Let $X$ be a geometrically connected variety of general type over a field $k$ finitely generated over $\Q$. Then $X(k)$ is not dense in $X$.
\end{BL}

Faltings proved that the weak Bombieri-Lang conjecture holds for sub-varieties of abelian varieties \cite[Theorem 5.1]{fal94}.

A variety $X$ over $k$ is \emph{geometrically mordellic} if every locally closed subvariety of $X_{\bar{k}}$ is of general type. If the weak Bombieri-Lang conjecture holds and $X$ is a geometrically mordellic variety over a field $k$ finitely generated over $\Q$, then it is immediate check that $X(k)$ is finite.

\begin{GL}
	Let $X$ be a smooth, projective, geometrically connected variety of general type over $\C$. There exists an open subset $U\s X$ which is geometrically mordellic.
\end{GL}

\begin{TA}
	Let $X$ be a smooth, hyperbolic curve over a field $K$ finitely generated over $\Q$. Assume that the section conjecture holds for all specializations of $X$ defined over number fields.
	
	Furthermore, assume either that
	\begin{itemize}
		\item the weak Bombieri-Lang conjecture holds over fields finitely generated over $\Q$, or
		\item the weak Bombieri-Lang conjecture holds over $\Q$ and the geometric Lang conjecture holds.
	\end{itemize}
	Then the section conjecture holds for $X$.
\end{TA}

Observe that Theorem A gives a complete splitting of these conjectures in statements either over number fields or over $\C$.

In \cite{bre21b} we proved that the weak Bombieri-Lang conjecture implies that if $f:X\to\P^{1}$ is a morphism over a field $k$ finitely generated over $\Q$ with geometrically mordellic fibers and $\P^{1}(k)\setminus f(X(k))$ is finite, then there exists a generic section $\spec k(\P^{1})\to X$. In \cite[Theorem C]{bre20b}, assuming the geometric Lang conjecture, we reduced the weak Bombieri-Lang conjecture to $\Q$. Theorem A follows from these results using known techniques in anabelian geometry.

\begin{remark*}
    Let us briefly mention a feature of the proof of Theorem A. There is some speculation about which varieties should be considered anabelian in dimension higher than $1$, since Grothendieck only gave examples and not definitions. In \autoref{gemcoh} we show that, conditional on Bombieri-Lang, geometrically mordellic varieties display anabelian behaviour: given a family of geometrically mordellic varieties, a Galois section of the generic fiber with geometric specializations is itself geometric. We only need this fact for curves, but the proof works in every dimension.
\end{remark*}

\subsection{New, unconditional results}

Without assuming the weak Bombieri-Lang conjecture, we prove the following.

\begin{TB}
	Let $X$ be a smooth, hyperbolic curve over a field $K$ finitely generated over $\Q$. Assume that the section conjecture holds for all specializations of $X$ defined over number fields. 
	
	Suppose that there exists a finite étale cover $Y\to X_{\bar{K}}$ and a non-trivial morphism $Y\to A_{\bar{K}}$, where $A$ is an abelian variety over $\bar{\Q}$. Then the section conjecture holds for $X$.
\end{TB}

An obvious question is the following: "how many" hyperbolic curves satisfy the hypothesis of Theorem B? We don't know any example of a curve which does not satisfy it. Let us give some examples of curves that do satisfy it.

\begin{CA}
	Fix an integer $g\ge 2$, and assume that the section conjecture holds for all projective, hyperelliptic curves of genus $g$ defined over number fields. Then it holds for all projective, hyperelliptic curves of genus $g$ defined over finitely generated extensions of $\Q$.
\end{CA}

\begin{CB}
	Let $X$ be a smooth, hyperbolic, affine curve over a field $K$ finitely generated over $\Q$. Assume that the section conjecture holds for all specializations of $X$ defined over number fields. Suppose that there exists a finite étale cover $Y\to X_{\bar{K}}$ and a non-trivial morphism $Y\to\P^{1}\setminus\{0,1,\infty\}$. Then the section conjecture holds for $X$. 
	
	In particular, if the section conjecture holds for all affine curves over all number fields, every curve defined over a finitely generated extension of $\Q$ contains an open subset which satisfies the section conjecture.
\end{CB}

\begin{CC}
	Fix integers $0\le g\le 2$ and $n\ge 0$ with $2g+n\ge 3$. Assume that the section conjecture holds for curves defined over number fields of the form $\bar{X}\setminus D$ where $\bar{X}$ is projective of genus $g$ and $D$ is a reduced divisor of degree $n$. Then the section conjecture holds for curves as above defined over finitely generated extensions of $\Q$. In particular, the section conjecture in genus $\le 2$ reduces to the case of number fields.
\end{CC}

\subsection{Discussion of contents}

We use abundantly well-known arguments in anabelian geometry, ideas that have been used by Nakamura, Tamagawa, Mochizuki, Hoshi, Sa\"idi, Pop, Stix among others. To name a few: using Deligne's weights to spread out Galois sections, using Hilbert's irreducibility theorem to show that certain Galois sections with geometric specializations are themselves geometric, Tamagawa's idea of neighbourhoods of a section and the reduction to curves with no rational points. We obviously do not claim originality for these arguments. To avoid any confusion, let us state clearly which are our original contributions.


The first is the idea is the connection between the weak Bombieri-Lang conjecture and the section conjecture. Secondly, thanks to a simplified approach avoiding Tate-Shafarevich groups, we observe that the techniques used by Sa\"idi in \cite{sai16} can prove a much stronger unconditional result. Third, we observe that this result applies to all hyperelliptic curves, all hyperbolic curves of genus $\le 2$ and an open subset of every curve.

In short, the main purpose of this article is not to introduce new techniques, but to prove new results using known techniques.

\subsection{Acknowledgements}

I would like to thank J. Stix for a discussion which helped me clarify the exposition of the paper and H. Esnault for many useful remarks. 

\subsection{Conventions}

We work over a base field $k$ of characteristic $0$. Curves are assumed to be smooth and geometrically connected over $k$. A gerbe over $k$ is finite if it has a finite groupoid presentation, or equivalently if its automorphism groups are finite.

We frequently use without mention the Mordell-Weil theorem about the rational points of abelian varieties and Faltings' theorem about the rational points of curves of genus $\ge 2$. These were originally proved for number fields, but they generalize to finitely generated extensions of $\Q$ \cite[Theorem 2.1]{con06} \cite[Theorem 5.1]{fal94}.

\section{Étale fundamental gerbes} We use the formalism of \emph{étale fundamental gerbes}, as developed by Borne and Vistoli \cite[\S 8]{bv15} \cite[Appendix A]{bre21}. Basically, it is an alternative point of view on the theory of étale fundamental groups. Given a geometrically connected scheme (or algebraic stack, or fibered category) $X$, the étale fundamental gerbe is a pro-finite étale gerbe $\Pi_{X/k}$ with a morphism $X\to\Pi_{X/k}$ which is universal with respect to morphisms to finite, étale gerbes over $k$. Its existence can be proved by a limit argument. Given a geometric point $x\in X(\bar{k})$ with image $s\in\Pi_{X/k}(\bar{k})$, there is a natural identification $\pi_{1}(X_{\bar{k}},x)=\aut_{\Pi_{X/k}}(x)$. If the image of $x$ is rational, then $\gal(\bar{k}/k)$ acts on $\aut_{\Pi_{X/k}}(x)$ and we have $\pi_{1}(X,x)=\aut_{\Pi_{X/k}}(x)\rtimes\gal(\bar{k}/k)$.

Loosely speaking, the étale fundamental gerbe is the space of Galois sections, but with a richer and more natural stacky structure with respect to the classical, set-theoretical one. In fact, if $X$ is geometrically connected over $k$ and $X\to\Pi_{X/k}$ is the étale fundamental gerbe, the space of Galois sections of $X$ is in natural bijection with the isomorphism classes of $\Pi_{X/k}(k)$, and more generally for a geometrically connected scheme $S$ over $k$ we have that $\Pi_{X/k}(S)/\sim$ corresponds to homomorphisms $\pi_{1}(S)\to\pi_{1}(X)$ over $\gal(\bar{k}/k)$ modulo the action of $\pi_{1}(X_{\bar{k}})$ by conjugation, see \cite[Proposition A.19]{bre21}. For instance, the structure morphism $X\to\Pi_{X/k}$ corresponds to the identity of the étale fundamental group of $X$. Given a rational point $\spec k\to X$, the associated Galois section $\spec k\to\Pi_{X/k}$ is obtained simply by composing with $X\to\Pi_{X/k}$. 

The two points of view are equivalent, and translations can be made back and forth. We use both, depending on which is convenient in each circumstance. Fundamental gerbes are naturally base point free and have simpler functorial properties, hence they are particularly well-suited for geometric constructions and specialization arguments in anabelian geometry. Fundamental groups are much more widely known, and most results in the literature are stated for groups, hence we use them when gerbes do not give particular advantages. The reader not familiar with fundamental gerbes might want to regard them just as spaces of Galois sections, and trust the fact that base points can be ignored.




In order to study family of curves using étale fundamental gerbes, we need to generalize the latter to a relative setting, i.e. define them not only over the spectrum of a field but over a more general basis. This will give us a global object containing all the spaces of Galois sections of the fibers in a coherent way. However, this is not the place for such a general task, hence we will define relative fundamental gerbes only for a small class of morphisms, exploiting the fact that smooth curves (except $\P^{1}$) have trivial second homotopy group.

\subsection{Good fibrations} 

Let $S$ be a smooth variety over $k$, $\bar{X}\to S$ a smooth projective morphism, $D\subset \bar{X}$ a divisor finite étale over $S$. Write $X=\bar{X}\setminus D$ and $F\s X_{\bar{k}}$ for a geometric fiber. Assume that $F,S$ are geometrically connected.

\begin{definition}
	A morphism $X\to S$ as above is a \emph{good fibration} over $k$ if the sequence
	\[1\to \pi_{1}(F)\to\pi_{1}(X_{\bar{k}})\to\pi_{1}(S_{\bar{k}})\to 1\]
	is exact.
\end{definition}


\begin{lemma}
	Let $\bar{X}\to S$, $D\s\bar{X}$, $X=\bar{X}\setminus D$ be as above. Assume that $S$ is a smooth curve, and $S_{\bar{k}}\not\simeq\P^{1}$. Then $X\to S$ is a good fibration
\end{lemma}

\begin{proof}
	With standard arguments, we may reduce to the case $k=\C$, so that étale fundamental groups are pro-finite completions of the topological ones. Since $\bar{X}\to S$ is projective, then $\pi_{1}^{\rm top}(F)$ is finitely generated. The short sequence of topological fundamental groups is exact since $\pi_{2}^{\rm top}(S)$ is trivial. The statement then follows from \cite[\S I.2.6, exercise 2.b]{ser94} and \cite[Proposition 3.6]{gjz08}.
\end{proof}

\subsection{Relative fundamental gerbes}

\begin{definition}\label{relgerbe}
	Define the étale fundamental gerbe $X\to\Pi_{X/S}$ of $X/S$ by the following $2$-cartesian diagram
	\[\begin{tikzcd}
		X\rar\ar[dr]\ar[rr,bend left]		&	\Pi_{X/S}\rar\dar\ar[dr,phantom,"\square",description]		&	\Pi_{X/k}\dar	\\
											&	S\rar														&	\Pi_{S/k}
	\end{tikzcd}\]
\end{definition}

\begin{lemma}\label{relgerbebc}
	Let $X\to S$ be a good fibration over $k$, and $S'$ a smooth, geometrically connected variety over some field extension $k'$ of $k$. Let $S'\to S$ be any morphism, and $X'\to S'$ the base change of $X\to S$. Then $X'\to S'$ is a good fibration over $k'$, and there is a natural isomorphism
	\[\Pi_{X'/S'}\xar{\sim}\Pi_{X/S}\times_{S}S'.\]
	\begin{proof}
		Choose the geometric fiber of $X\to S$ in the image of $X'\to X$. The fact that $X'\to S'$ is a good fibration is straightforward. The étale fundamental gerbe behaves well with respect to base change  \cite[Proposition A.23]{bre21}, hence $\Pi_{S\times k'/k'}\simeq\Pi_{S/k}\times k'$ and similarly for $X$. There is a $2$-commutative diagram
		\[\begin{tikzcd}[column sep=tiny,row sep=tiny]
																	&	\Pi_{X\times k'/S\times k'}\ar[rr]\ar[dd]		&										&	\Pi_{X\times k'/k'}\ar[dd]	\\
			\Pi_{X'/S'}\ar[rr,crossing over]\ar[ur,dotted]\ar[dd]	&								&	\Pi_{X'/k'}\ar[ur]					&						\\
																	&	S\times k'\ar[rr]					&										&	\Pi_{S\times k'/k'}			\\
			S'\ar[ur]\ar[rr]										&				&	\Pi_{S'/k'}\ar[ur]\ar[leftarrow,uu,crossing over]	&
		\end{tikzcd}\]
		where the front and back squares are $2$-cartesian by definition. Since there is a commutative diagram with exact rows
		\[\begin{tikzcd}
			1\rar	&	\pi_{1}(F)\rar\dar[equal]	&	\pi_{1}(X'_{\bar{k'}})\rar\dar	&	\pi_{1}(S'_{\bar{k'}})\rar\dar	& 1	\\
			1\rar	&	\pi_{1}(F)\rar				&	\pi_{1}(X_{\bar{k}})\rar		&	\pi_{1}(S_{\bar{k}})\rar		& 1
		\end{tikzcd}\]
		then $\pi_{1}(X'_{\bar{k'}})\simeq \pi_{1}(X_{\bar{k'}})\times_{\pi_{1}(S_{\bar{k'}})}\pi_{1}(S'_{\bar{k'}})$ and hence $B\pi_{1}(X'_{\bar{k'}})\simeq B\pi_{1}(X_{\bar{k'}})\times_{B\pi_{1}(S_{\bar{k'}})}B\pi_{1}(S'_{\bar{k'}})$, where $BG$ is the classifying stack of $G$. The base change of $\Pi_{S/k}$ to $\bar{k'}$ is isomorphic to $B\pi_{1}(S_{\bar{k'}})$, and similarly for $X,S',X'$. This implies that the right square is $2$-cartesian: by descent theory, we can check it after base changing to $\bar{k'}$. By diagram chasing we get that the dotted arrow $\Pi_{X'/S'}\to\Pi_{X\times k'/S\times k'}$ exists and that the left square is $2$-cartesian. This completes the proof.
	\end{proof}
\end{lemma}

\section{Spreading out sections}\label{sprout}

Let $C$ be an affine curve over $k$ and $X\to C$ a good fibration, we study when a generic section $\spec k(C)\to\Pi_{X_{k(C)/k(C)}}$ spreads out to a global section $C\to\Pi_{X/C}$, or equivalently when an homomorphism $\gal(\overline{k(C)}/k(C))\to\pi_{1}(X)$ over $\pi_{1}(C)$ extends to a section $\pi_{1}(C)\to \pi_{1}(X)$. This is a well-known topic, see \cite[Chapter 8]{sti13} \cite[\S 2]{sai16}. We can't find a precise reference for the statements we need, though, hence let us give proofs. The arguments are well-known, we do not claim originality.

In this section, the point of view of gerbes does not simplify the arguments, hence in the proofs we use groups for the convenience of the reader. If the reader prefers the point of view of gerbes, \cite[Appendix A]{bre21c} contains the relevant techniques.

\begin{definition}
	Let $G$ be a pro-finite group. We say that $G$ \emph{has no abelian torsion} if, for every open subgroup $H\s G$, the abelianization of $H$ is torsion free.
\end{definition}

\begin{example}
	The étale fundamental groups of smooth, projective curves and abelian varieties over algebraically closed fields have no abelian torsion. For abelian varieties, this follows from the fact that their fundamental group is abelian and torsion free. For a smooth, projective curve, it follows from the fact that its abelianized fundamental group coincides with the fundamental group of the Jacobian \cite[Proposition 69]{sti13}, plus the fact that a finite étale cover of a smooth, projective curve is again a smooth, projective curve.
\end{example}

\begin{lemma}\label{trivab}
	Let $k$ be a field finitely generated over $\Q$, and $X$ a smooth, projective, geometrically connected variety such that $\pi_{1}(X_{\bar{k}})$ has no abelian torsion. Assume furthermore that $X(k)\neq\emptyset$ and choose a rational base point, so that $\gal(\bar{k}/k)$ acts on $\pi_{1}(X_{\bar{k}})$. Every Galois-equivariant homomorphism $\hz(1)\to\pi_{1}(X_{\bar{k}})$ is trivial.
\end{lemma}

\begin{proof}
	It is enough to show that every Galois-equivariant homomorphism $\Z_{\ell}(1)\to\pi_{1}(X_{\bar{k}})$ is trivial for every prime $\ell$. Assume first that $X=A$ is an abelian variety, we want to show that a Galois-equivariant homomorphism $\Z_{\ell}(1)\to T_{\ell}A$ is trivial. Since $k$ is finitely generated over $\Q$, we may find a scheme $V$ of finite type over $\Z$ and an abelian scheme $\mc{A}\to V$ whose generic fiber is $A$. Then $T_{\ell}A$ has weight $-1$ since it is dual to $\H^{1}(A_{\bar{k}},\Z_{\ell})$ while $\Z_{\ell}(1)$ has weight $-2$, hence $\Z_{\ell}(1)\to T_{\ell}A$ is trivial.
	
	Let us now do the general case. If, by contradiction, $\hz(1)\to\pi_{1}(X_{\bar{k}})$ is not trivial, then up to passing to a finite extension of $k$ and a finite, étale covering of $X$ we may assume that the composition $\hz(1)\to\pi_{1}(X_{\bar{k}})\to\pi_{1}^{\rm ab}(X_{\bar{k}})$ is non-trivial. By \cite[Proposition 69]{sti13}, $\pi_{1}^{\rm ab}(X_{\bar{k}})=\pi_{1}(A_{X,\bar{k}})$, where $A_{X}$ is the Albanese variety of $X$. We conclude by the preceding case. 
\end{proof}

\begin{corollary}\label{curprop}
	Let $C$ an affine curve over a field $k$ finitely generated over $\Q$ and $X\to C$ a smooth, projective morphism with geometrically connected fibers. Suppose that the étale fundamental group of the fibers has no abelian torsion. Every generic section $\spec k(C)\to\Pi_{X_{k(C)}/k(C)}$ extends to a unique section $C\to\Pi_{X/C}$.
\end{corollary}

\begin{proof}
	Equivalently, we have to show that for every homomorphism $\gal(\overline{k(C)}/k(C))\to\pi_{1}(X)$ over $\pi_{1}(C)$, every finite extension $h/k$ and every $c\in C(h)$, the inertia subgroup $\hz(1)\subset\gal(\overline{k(C)}/k(C))$ maps to $0$ in $\pi_{1}(X)$.
	
	Up to a finite extension of $k$, we may assume $h=k$ and $X_{c}(k)\neq\emptyset$, hence we may write $\pi_{1}(X_{c})\simeq\pi_{1}(X_{c,\bar{k}})\rtimes\gal(\bar{k}/k)$. Since the inertia maps to zero in $\pi_{1}(C)$, we have a Galois-equivariant factorization $\hz(1)\to\pi_{1}(X_{c,\bar{k}})\to\pi_{1}(X_{\bar{k}})$, which is trivial by \autoref{trivab}.
\end{proof}

\section{The coherence property}

\autoref{curprop} allows us to \emph{specialize} Galois sections for families of projective curves or abelian varieties. In view of this, the passage from number fields to finitely generated fields for projective curves essentially reduces to the following \emph{coherence property}.

\begin{definition}\label{cohdef}
	Let $C$ be an affine curve over $k$ and $X\to C$ a good fibration. A \emph{quasi-section} of $\Pi_{X/C}$ is the datum of an affine curve $C'$, a quasi-finite morphism $C'\to C$ and a section $C'\to\Pi_{X/C}$. A quasi-section $s:C'\to\Pi_{X/C}$ is \emph{pointwise geometric} if, for every finite extension $h/k$ and every point $c\in C'(h)$, the restriction of $s$ to a section $s_{c}\in\Pi_{X_{c}/h}(h)$ is geometric.
	
	We say that $X\to C$ has the \emph{coherence property} if, for every pointwise geometric quasi-section $s:C'\to\Pi_{X/C}$, the restriction to the generic point $s_{k(C')}$ is geometric.
\end{definition}


\subsection{Fibrations of geometrically mordellic varieties}


If the weak Bombieri-Lang conjecture holds, a higher dimensional generalization of Hilbert's irreducibility theorem is true, and this implies that fibrations of geometrically mordellic varieties have the coherence property.

\begin{theorem}[{\cite[Theorem B]{bre21b}}]\label{khit}
	Assume that the weak Bombieri-Lang conjecture holds, and let $k$ be a field finitely generated over $\Q$. Let $f:X\to C$ be a morphism of varieties over $k$, with $C$ a geometrically connected curve, and assume that the fibers of $f$ are geometrically mordellic. If $X(h)\to C(h)$ is surjective for every finite extension $h/k$, there exists a generic section $\spec k(C)\to X$. 
\end{theorem}


\begin{proposition}\label{gemcoh}
	Let $C$ be an affine curve and $X\to C$ a good fibration over a field $k$ finitely generated over $\Q$. Assume that the fibers of $X\to C$ are geometrically mordellic and that the weak Bombieri-Lang conjecture holds. Then $X\to C$ has the coherence property.
	\begin{proof}
		Let $s:C'\to\Pi_{X/C}$ be a pointwise geometric quasi-section, we want to prove that $s_{k(C')}$ is geometric. By base change, may assume $C'=C$. Since the fibers of $X\to C$ are geometrically mordellic and we are assuming that the weak Bombieri-Lang conjecture holds, it follows that $X(k(C))$ is finite. Proving that $s_{k(C)}\in\Pi_{X/C}(k(C))$ is geometric is then equivalent to proving that the composition $\spec k(C)\to\Pi_{X/C}\to\Pi_{X/k}$ is associated with one of the finite number of sections $\spec k(C)\to X$.
		
		Assume by contradiction that this is not true, then there exists a morphism $\Pi_{X/k}\to\Psi$ with $\Psi$ a finite étale gerbe such that the composition $C\to\Psi$ is not associate with any of the sections $\spec k(C)\to X$. Equivalently, the fibered product $X'=X\times_{\Psi\times C}C$ satisfies $X'(k(C))=\emptyset$. 
		
		Let $h/k$ be a finite extension and $c\in C(h)$ a point, since $s$ is pointwise geometric the specialization $s_{c}\in\Pi_{X_{c}/h}(h)$ is associated with a point $x\in X_{c}(h)$. It follows that $(x,c)\in X\times_{\Psi\times C}C(h)$ defines a point of $X'(h)$ over $c\in C(h)$. In particular, $X'(h)\to C(h)$ is surjective for every finite extension $h/k$. This is in contradiction with \autoref{khit}, since the fact that the fibers of $X\to C$ are geometrically mordellic implies the same for $X'\to C$.
	\end{proof}
\end{proposition}

\subsection{Extension of scalars}

The following is a well-known lemma (see for instance \cite[Proposition 109]{sti13}), but we can't find a reference for the exact statement that we need.

\begin{lemma}\label{finite}
	Let $k$ be a subfield of a finite extension of $\Q_{p}$, and $X$ either an abelian variety or a (possibly non-projective) curve. Let $k'/k$ be a finite extension, and $s\in\Pi_{X/k}(k)$ a Galois section such that $s_{k'}$ is geometric (resp. cuspidal). Then $s$ is geometric (resp. cuspidal).
\end{lemma}

\begin{proof}
	Let us show that every Galois section of $X$ has trivial centraliser cf. \cite[Definition 11]{sti13}. By definition of centraliser, it is enough to prove this over some extension of $k$, we may thus assume that $k$ is a finite extension of $\Q_{p}$. If $X$ is a curve, then this is \cite[Lemma 104]{sti13}. If $X$ is an abelian variety, then $\pi_{1}(X)=\gal(\bar{k}/k)\rtimes TX$ and centralisers of sections are in bijection with $\H^{0}(k,TX)$, which is trivial by \cite[Lemma 105]{sti13}. Triviality of centralisers implies that the space of Galois sections of $X$ as a functor from finite extensions of $k$ to sets is a sheaf for the étale topology \cite[Proposition 28]{sti13}. It is then enough to prove that geometric and cuspidal sections are subsheaves.
	
	Since the étale topology is subcanonical, the fact that geometric sections are a subsheaf follows from the injectivity of the map $X(k)\to\Pi_{X/k}(k)$ \cite[Proposition 75]{sti13}. This completes the proof in the projective case. Assume that $X$ is a curve with completion $\bar{X}$.
	
	The set of cuspidal sections has the form $\sqcup_{x\in\bar{X}\setminus X(k)}\H^{1}(k,\hz(1))$, see for instance \cite[\S 8]{bre21}. Hence, it is enough to show that $\H^{1}(k,\hz(1))$ is a sheaf. By the Hochschild-Serre spectral sequence (or by the identification $\H^{1}(k,\hz(1))$ with the rational points of the classifying stack $B\hz(1)$) this reduces to showing that $\H^{0}(k,\hz(1))=\projlim_{n}\mu_{n}(k)=0$, which is true if $k$ is a finite extension of $\Q_{p}$ (and hence for its subfields).
\end{proof}

\subsection{Fibrations of abelian varieties}\label{sect:cohab}


Given an abelian group $G$, write $\hat{G}=\projlim_{n} G/nG$.

\begin{lemma}\label{hatcar}
	Let $G\s H$ finitely generated abelian groups. Then $G=\hat{G}\cap H\s \hat{H}$, or equivalently $\hat{G}/G\to\hat{H}/H$ is injective.
	\begin{proof}
		This follows from the fact that $\hz/\Z$ is a flat $\Z$-module and $\hat{G}=G\otimes_{\Z}(\hz/\Z)$.
			
	\end{proof}
\end{lemma}

If $A$ is an abelian variety over $k$, for every positive integer $n$ Kummer theory gives us an injective map $A(k)/nA(k)\hookrightarrow\H^{1}(k,A[n])$, passing to the limit we obtain an injective map $\hat{A(k)}\hookrightarrow\H^{1}(k,TA)$. The space of Galois sections of $A$ is in natural bijection with $\H^{1}(k,TA)$, and the map sending a rational point to its associated Galois section coincides with $A(k)\to\hat{A(k)}\to\H^{1}(k,TA)$ \cite[Proposition 71]{sti13}. If $k$ is finitely generated over $\Q$, by the Mordell-Weil theorem the homomorphism $A(k)\to\hat{A(k)}$ is injective, hence we have embeddings $A(k)\s \hat{A(k)}\s \H^{1}(k,A)=\Pi_{A/k}(k)$.

\begin{proposition}\label{cohab}
	Let $k$ be a field finitely generated over $\Q$. Let $C$ be an affine curve and $A$ an abelian variety over $k$. Then $A\times C\to C$ has the coherence property.
	\begin{proof}
		Suppose that we have a pointwise geometric quasi-section $C'\to\Pi_{A\times C/C}$, we want to show that $s_{k(C')}$ is geometric. By base change we may assume that $C'=C$. Furthermore, since $\Pi_{A\times C/C}=\Pi_{A/k}\times C$, we have just a section $s:C\to\Pi_{A/k}$ with geometric specializations. By \autoref{finite}, it is enough to prove that $s_{k(C)}$ is geometric after a finite extension of $k$.
		
		Let $J$ be the Jacobian variety of $C$, since $C$ generates $J$ up to extending $k$ we may assume that there exists a subgroup $G\s J(k)$ generated by a finite number of rational points of $C$ and which is dense in $J$. In particular,
		\[\hom(J,A)\hookrightarrow\hom(G,A(k))\]
		is injective.
		
		Let $J_{0}$ be the semi-abelian Jacobian of $C$, its Tate module is the abelianization of the étale fundamental group of $C$, hence we have an induced homomorphism $TJ_{0}\to TA$. Denote by $J_{t}\s J_{0}$ the toric part, by weight reasons the restriction $TJ_{t}\to TA$ is trivial, hence there is an induced homomorphism $\phi:TJ\to TA$.
		
		By Faltings' theorem, this homomorphism is associated with an element of $\widehat{\hom(J,A)}$.	Let us show that $\phi$ is associated with an actual morphism $J\to A$. By \autoref{hatcar},
		\[\hom(J,A)=\hom(G,A(k))\cap \widehat{\hom(J,A)}\s \widehat{\hom(G,A(k))}=\hom(G,\hat{A(k)}),\]
		hence it is enough to show that $\psi(G)\subset A(k)\subset \hat{A(k)}$, where $\psi\in\hom(G,\hat{A(k)})$ is the homomorphism associated with $\phi$.{}
			
		Let $g_1,\dots,g_n\in C(k)$ be a finite set of generators for $G$. By \cite[Proposition 71]{sti13}, since $g_{i}\in C(k)$, the image of $\psi(g_{i})$ in $\H^{1}(k,TA)=\Pi_{A/k}(k)$ coincides with the specialization $s_{g_{i}}$. Since $s_{g_{i}}$ is geometric by hypothesis, $\psi(g_{i})\in A(k)$ for every $i$, hence $\psi(G)\s A(k)$.
	\end{proof}
\end{proposition}

\subsection{Pulling the coherence property along a finite morphism}


The following is a well-known consequence of Hilbert's irreducibility theorem, but we could not find a suitable reference.

\begin{lemma}\label{hilb}
	Let $k$ be an Hilbertian field of characteristic $0$, and $f:Y\to X$ a quasi-finite morphism of schemes of finite type over $k$ of dimension $1$, with $X$ integral. Suppose that $Y(h)\to X(h)$ is surjective for every finite extension $h/k$. Then there is a generic section $\spec k(X)\to Y$.
	\begin{proof}
		Up to shrinking $X$ and $Y$, we may assume that $Y\to X$ is finite and that there exists a finite, flat morphism $X\to U$ where $U$ is an open subset of $\A^{1}$. The Weil restriction $R_{X/U}Y$ exists \cite[Theorem 7.6.4]{blr90}, is finite over $U$ and the hypothesis implies that $R_{X/U}Y(k)\to U(k)$ is surjective. The fact that $Y\to X$ is finite implies that $R_{X/U}Y\to U$ is finite, too (in general, the relative dimension of $R_{X/U}Y/U$ is the relative dimension of $Y/X$ multiplied by $\deg X/U$). The hypothesis implies that $R_{X/U}Y(k)\to U(k)$ is surjective. By Hilbert's irreducibility theorem there exists a section $\spec k(U)\to R_{X/U}Y$ which induces a section $\spec k(X)\to Y$.
	\end{proof}
\end{lemma}

\begin{lemma}\label{geogalsep}
	Let $C$ be an affine curve over a field $k$ finitely generated over $\Q$, $\pi:X\to C$ a family of smooth, proper curves of genus at least $2$. Let $Y\to C$ a good fibration with a finite morphism $\omega:X\to Y$ over $C$.
	
	If $Y\to C$ has the coherence property, then $X\to C$ has the coherence property too.
	\begin{proof}
		Let $s:C'\to\Pi_{\mc{X}/C}$ be a pointwise geometric quasi-section, we want to show that $s_{k(C')}$ is geometric. By base change, we may assume $C'=C$. Since $Y\to C$ has the coherence property, up to shrinking $C$ we have that $\omega\circ s$ is associated with a section $f:C\to Y$.	Assume by contradiction that $s_{k(C)}$ is not geometric, we may then find a morphism $\Pi_{X/k}\to\Psi$ with $\Psi$ is a finite étale gerbe such that $X'=X\times_{\Psi\times C}C$ satisfies $X'(k(C))=\emptyset$ as in the proof of \autoref{gemcoh}.
		
		Let $D\s X$, $D'\s X'$ be the inverse images of $f(C)$. Since $X'\to X\to Y$ is finite, $D'\to f(C)=C$ is finite. Let $h/k$ be a finite extension and $c\in C(h)$ a point, then by hypothesis $s_{c}\in\Pi_{X_{c}/h}$ is associated with a point $x\in X_{c}(h)$, and clearly $x\in D(h)$. The pair $(x,c)\in X\times_{\Psi\times C}C(h)=X'(h)$ defines a point of $D'(h)$ over $C(h)$, hence $D'(h)\to C(h)$ is surjective. By \autoref{hilb}, $D'(k(C))\s X'(k(C))$ is non-empty, which is absurd.
	\end{proof}
\end{lemma}

\begin{corollary}\label{abvred}
	Let $\pi:X\to C$ a family of smooth, proper curves of genus at least $2$, with $C$ an affine curve over a field $k$ finitely generated over $\Q$. If there exists an abelian variety $A$ over $k$ with a non-trivial morphism $X_{k(C)}\to A_{k(C)}$ then $X\to C$ has the coherence property.\qed
\end{corollary}

\section{Proofs of the main theorems}

Up to minor differences, the proof of Theorem A is strictly contained in that of Theorem B. Because of this, let us prove Theorem B first.

\begin{TB}
	Let $X$ be a smooth, hyperbolic (not necessarily projective) curve over a field $K$ finitely generated over $\Q$. Assume that the section conjecture holds for all specializations of $X$ defined over number fields. Suppose that there exists a finite étale cover $Y\to X_{\bar{K}}$ and a non-trivial morphism $Y\to A_{\bar{K}}$, where $A$ is an abelian variety over $\bar{\Q}$. Then the section conjecture holds for $X$.
\end{TB}

\begin{proof}
	By induction, we may assume that the theorem holds for finitely generated extensions of $\Q$ of transcendence degree strictly smaller than $\on{trdeg}(K/\Q)$. Let $s\in\Pi_{X/K}(K)$ be a section, assume by contradiction that $s$ does not come from geometry. By \autoref{finite}, we may replace $K$ with an arbitrary finite extension.
	
	By hypothesis, there exists a finite extension $K'/K$, a finite étale cover $Y\to X_{K'}$ and a non-trivial morphism $Y\to A_{K'}$, where $A$ is an abelian variety over a number field contained in $K'$, and we may assume that $Y$ has genus $\ge 2$. Up to a further finite extension of $K'$, we may assume that $s_{K'}$ lifts to a section $s'\in\Pi_{Y/K'}(K')$. Standard arguments about finite étale coverings and scalar extensions in anabelian geometry \cite[Propositions 109, 110, 111]{sti13} imply that the section conjecture holds for the specializations of $Y$ defined over number fields, too. By \autoref{finite}, $s_{K'}$ and thus $s'$ do not come from geometry, hence we may replace $K,X$ with $K',Y$ and assume that there exists a non-trivial morphism $X\to A_{K}$.
	
	Furthermore, up to replacing $X$ with an étale neighbourhood of $s$, we may assume that the completion $\bar{X}$ of $X$ has no rational points, see the proof of \cite[Proposition 103]{sti13}. The reference works over a number field in order to use the Mordell-Weil theorem for abelian varieties, which guarantees the injectivity of the pro-finite Kummer map, and Faltings' theorem for curves of genus $\ge 2$, which is used in order to apply an argument by Tamagawa \cite[Corollary 99]{sti13}. These theorems both generalize to finitely generated extensions of $\Q$ \cite[Theorem 2.1]{con06} \cite[Theorem 5.1]{fal94}, hence the argument generalizes word by word to finitely generated extensions of $\Q$.
		
	Let $k\s K$ be a subfield algebraically closed in $K$ and such that $\on{trdeg}(K/k)=1$, and $C$ an affine curve over $k$ such that $k(C)=K$. Up to shrinking $C$, we may assume that $\bar{X}\to \spec K$ spreads out to a family of smooth, proper hyperbolic curves $\bar{\mc{X}}\to C$. Rational maps from smooth varieties to abelian varieties extend to morphisms, hence there is an extension $\bar{\mc{X}}\to A_{k}$. Up to shrinking $C$, we may assume that there exists a divisor $D\subset \bar{\mc{X}}$ finite étale over $C$ such that the generic fiber of $\mc{X}=\bar{\mc{X}}\setminus D$ is $X$.
	
	Let $\bar{s}\in\Pi_{\bar{X}/K}(K)$ be the image of $s$, by \autoref{curprop} it spreads out to a section $\tilde{s}\in\Pi_{\bar{\mc{X}}/C}(C)$. By \cite[\S 2]{bre21d}, we can define (non-unique) specializations of $s$, and these map to specializations of $\tilde{s}$. Since specialization of curves is transitive, the closed fibers of $\mc{X}\to C$ satisfy the hypothesis of the theorem, hence by induction hypothesis the section conjecture holds for them. In particular, the specializations of $s$ come from geometry, hence the specializations of $\tilde{s}$ are geometric. By \autoref{abvred} $\bar{\mc{X}}\to C$ has the coherence property and hence $\bar{s}_{K}$ is geometric, which is absurd since $\bar{X}$ has no rational points.
\end{proof}

\begin{CA}
	Fix an integer $g\ge 2$, and assume that the section conjecture holds for all projective, hyperelliptic curves of genus $g$ defined over number fields. Then it holds for all projective, hyperelliptic curves of genus $g$ defined over finitely generated extensions of $\Q$.
\end{CA}

\begin{proof}
	By \cite[Proposition 1.8]{bt02}, every hyperelliptic curve has a finite étale covering which maps to a curve of genus $2$ defined over $\bar{\Q}$. By looking at the Jacobian of this curve, we see that every hyperelliptic curve satisfies the hypothesis Theorem B.
\end{proof}

\begin{CB}
	Let $X$ be a smooth, hyperbolic, affine curve over a field $K$ finitely generated over $\Q$. Assume that the section conjecture holds for all specializations of $X$ defined over number fields. Suppose that there exists a finite étale cover $Y\to X_{\bar{K}}$ and a non-trivial morphism $Y\to\P^{1}\setminus\{0,1,\infty\}$. Then the section conjecture holds for $X$. 
	
	In particular, if the section conjecture holds for all affine curves over all number fields, every curve defined over a finitely generated extension of $\Q$ contains an open subset which satisfies the section conjecture.
\end{CB}

\begin{proof}
	The elliptic curve $E$ over $\bar{\Q}$ of $j$ invariant $0$ has a triple cover $E\to\P^{1}$ branched only over $0,1,\infty$, hence $X$ satisfies the hypothesis of Theorem B.
\end{proof}

\begin{CC}
	Fix integers $0\le g\le 2$ and $n\ge 0$ with $2g+n\ge 3$. Assume that the section conjecture holds for curves defined over number fields of the form $\bar{X}\setminus D$ where $\bar{X}$ is projective of genus $g$ and $D$ is a reduced divisor of degree $n$. Then the section conjecture holds for curves as above defined over finitely generated extensions of $\Q$. In particular, the section conjecture in genus $\le 2$ reduces to the case of number fields.
\end{CC}

\begin{proof}
	If $g=0$, we may clearly apply Corollary B.
	
	If $E$ is an elliptic curve, there exist finite étale morphisms $E\setminus E[2]\to E\setminus\{0\}$ and $E\setminus E[2]\to\P^{1}\setminus\{0,1,\infty,\lambda\}$, where $E[2]$ is the $2$-torsion and $\lambda$ is some constant. Because of this, the case $g=1$ follows from Corollary B.
	
	If $g=2$, then we have shown in the proof of Corollary A that every curve of genus $2$ satisfies the hypothesis of Theorem B. 
\end{proof}

\begin{TA}
	Let $X$ be a smooth, hyperbolic curve over a field $K$ finitely generated over $\Q$. Assume that the section conjecture holds for all specializations of $X$ defined over number fields.
	
	Furthermore, assume either that
	\begin{itemize}
		\item the weak Bombieri-Lang conjecture holds over fields finitely generated over $\Q$, or
		\item the weak Bombieri-Lang conjecture holds over $\Q$ and the geometric Lang conjecture holds.
	\end{itemize}
	Then the section conjecture holds for $X$.
\end{TA}

\begin{proof}
	By \cite[Theorem C]{bre20b}, the second condition implies the first, we may thus assume that the first holds. The proof is now analogous to that of Theorem B, ignoring the part about abelian varieties and replacing \autoref{abvred} with \autoref{gemcoh}.

	
\end{proof}

\printbibliography
	
\end{document}